\documentclass[11pt,reqno,tbtags]{amsart}
\pdfoutput=1
\usepackage[]{amsmath,amssymb,amsfonts,latexsym,amsthm,enumerate}
\usepackage{amssymb,bbm,mathrsfs}
\usepackage{enumerate,graphicx,paralist}
\usepackage[numeric,initials,nobysame]{amsrefs}

\numberwithin{equation}{section}
\usepackage[colorlinks=true, pdfstartview=FitV, linkcolor=blue, citecolor=blue, urlcolor=blue,pagebackref=false]{hyperref}
\usepackage[marginratio=1:1,height=584pt,width=488pt,tmargin=117pt]{geometry}
\usepackage[margin=25pt,font=small,justification=centering
]{caption}

\newtheorem{maintheorem}{Theorem}

\newtheorem{theorem}{Theorem}[section]
\newtheorem*{theorem*}{Theorem}
\newtheorem{lemma}[theorem]{Lemma}

\newtheorem{corollary}[theorem]{Corollary}

\theoremstyle{definition}{

\newtheorem*{definition*}{Definition}

\newtheorem*{example*}{Example}

\newtheorem*{remark*}{Remark}
}

\newcommand{\Z}{\mathbb Z}

\newcommand{\E}{\mathbb{E}}
\renewcommand{\P}{\mathbb{P}}
\DeclareMathOperator{\var}{Var}
\DeclareMathOperator{\Cov}{Cov}

\newcommand{\tmix}{t_\textsc{mix}}

\newcommand{\tv}{{\textsc{tv}}}

\newcommand{\given}{\, \big| \,}

\newcommand{\one}{\mathbbm{1}}
\newcommand{\red}{\textsc{Red}}
\newcommand{\blue}{\textsc{Blue}}
\newcommand{\green}{\textsc{Green}}

\renewcommand{\epsilon}{\varepsilon}
\renewcommand{\phi}{\varphi}

\newcommand{\cC}{\mathcal{C}}

\newcommand{\cF}{\mathcal{F}}

\newcommand{\cU}{\mathcal{U}}

\newcommand{\cX}{\mathcal{X}}

\newcommand{\len}{{\mathfrak{L}}}
\newcommand{\sm}{{\mathfrak{m}}}
\newcommand{\anim}{{\mathfrak{W}}}

\newcommand{\sH}{\mathscr{H}}

\newcommand{\tcut}{t_\sm}
\newcommand{\scut}{s_\star}
\newcommand{\tpluss}{t_\star}
\newcommand{\tminuss}{t_\star^-}

\date{}

\begin{document}
\title{An exposition to information percolation for the Ising model}

\author{Eyal Lubetzky}
\address{Eyal Lubetzky\hfill\break
Courant Institute 
\\ New York University\\
251 Mercer Street\\ New York, NY 10012, USA.}
\email{eyal@courant.nyu.edu}
\urladdr{}

\author{Allan Sly}
\address{Allan Sly\hfill\break
Department of Statistics\\
UC Berkeley\\
Berkeley, CA 94720, USA.}
\email{sly@stat.berkeley.edu}
\urladdr{}

\begin{abstract}
Information percolation is a new method for analyzing stochastic spin systems through classifying and controlling the clusters of information-flow in the space-time slab. It yielded sharp mixing estimates (cutoff with an $O(1)$-window) for the Ising model on $\Z^d$ up to the critical temperature, 
as well as results on the effect of initial conditions on mixing.
In this expository note we demonstrate the method on lattices (more generally, on any locally-finite transitive graph) at very high temperatures.
\end{abstract}
\maketitle
\vspace{-0.75cm}

\section{Introduction}

The Ising model on a finite graph $G$ with vertex-set $V$ and edge-set $E$ is a distribution over the set of configurations $\Omega=\{\pm1\}^V$;
each $\sigma\in\Omega$ is an assignment of plus/minus \emph{spins} to the sites in $V$, and the probability of $\sigma \in \Omega$ is given by the Gibbs distribution
\begin{equation}
  \label{eq-Ising}
  \pi(\sigma)  = \mathcal{Z}^{-1} e^{\beta \sum_{uv\in E} \sigma(u)\sigma(v) } \,,
\end{equation}
where $\mathcal{Z}$ is a normalizer (the partition-function) and $\beta$ is the inverse-temperature, here taken to be non-negative (ferromagnetic). These definitions extend to infinite locally finite graphs (see, e.g.,~\cites{Liggett,Martinelli97}).
The (continuous-time) heat-bath Glauber dynamics for the Ising model is the Markov chain---reversible w.r.t.\ the Ising measure $\pi$---where each site is associated with a rate-1 Poisson clock, and as the clock at some site $u$ rings, the spin of $u$ is replaced by a sample from the marginal of $\pi$ given all other spins.

An important notion of measuring the convergence of a Markov chain $(X_t)$ to its stationarity measure $\pi$ is its total-variation mixing time, denoted $\tmix(\epsilon)$ for a precision parameter $0<\epsilon<1$:
\[ \tmix(\epsilon) = \inf\big\{t \;:\; \max_{x_0 \in \Omega} \| \P_{x_0}(X_t \in \cdot)- \pi\|_\tv \leq \epsilon \big\}\,,\]
where here and in what follows $\P_{x_0}$ denotes the probability given $X_0=x_0$, and the total-variation distance
$
\|\nu_1-\nu_2\|_\tv$ is defined as $\max_{A\subset \Omega} |\nu_1(A)-\nu_2(A)| = \tfrac12\sum_{\sigma\in\Omega} |\nu_1(\sigma)-\nu_2(\sigma)|$.

The impact of the precision parameter $\epsilon$ in this definition is addressed by the {\em cutoff phenomenon}---a concept going back to the pioneering works~\cites{Aldous,AD,DiSh}---roughly saying the choice of any fixed $\epsilon$ does not change the asymptotics of $\tmix(\epsilon)$ as the system size goes to infinity.
Formally, a family of ergodic finite Markov chains $(X_t)$, indexed by an implicit parameter $n$, exhibits \emph{cutoff} iff $\tmix(\epsilon)=(1+o(1))\tmix(\epsilon')$ for any fixed $0<\epsilon,\epsilon'<1$. The \emph{cutoff window} addresses the correction terms: a sequence $w_n = o\big(\tmix(1/2)\big)$ is a cutoff window if $\tmix(\epsilon) = \tmix(1-\epsilon) + O(w_n)$ for any $0<\epsilon<1$ with an implicit constant that may depend on $\epsilon$. That is, the Markov chain exhibits a sharp transition in its convergence to equilibrium, whereby its distance drops abruptly (along the cutoff window) from near 1 to near 0.

Establishing cutoff can be highly challenging (see the survey~\cite{Diaconis}), even for simple random walk with a uniform stationary measure: e.g., it is conjectured that on \emph{every} transitive expander graph the random walk exhibits cutoff, yet there is not a single example of a transitive expander where cutoff was confirmed (even without transitivity, there were no examples of expanders with cutoff before~\cites{LS2,LSexp}). As for Glauber dynamics for the Ising model, where our understanding of $\pi$ is far more limited, until recently cutoff was confirmed only in a few cases, with first results on lattices appearing in the works~\cites{LS1,LS3}.

The methods used to analyze the dynamics for the Ising model on $\Z^d$ in those works had several caveats: the reliance on delicate features such as log-Sobolev inequalities did not cover the full high temperature regime in dimensions $d\geq 3$, and did not give a correct bound on the cutoff window; furthermore, the argument immediately broke on any geometry with exponential ball growth (expanders). Recently, in~\cite{LS4} and its companion paper~\cite{LS5}, we introduced a new framework called \emph{information percolation}, which does not have these limitations and can hopefully be used to analyze a wide range of stochastic spin systems. We demonstrated its application in analyzing the stochastic Ising model on $\Z^d$ up to $\beta_c$, and to compare the effects of different starting configurations on mixing---e.g., showing that an initial state of i.i.d.\ spins mixes about twice faster than the all-plus starting state (which is essentially the worst one), while almost every starting state is as bad as the worst one.

Here we demonstrate a simpler application of the method for the Ising model on a fixed-degree transitive graph at very high temperatures, establishing cutoff within an $O(1)$-window around
\begin{equation}\label{eq-t*-def}
\tcut = \inf\left\{\; t>0 \;:\;  \sm_t \leq 1/\sqrt{n} \;\right\}\,,
\end{equation}
where $\sm_t = \E X_t^+(v)$ is the magnetization at the origin at time $t>0$ (in which $X_t^+$ denotes the dynamics started from all-plus). That is, $\tcut$ is the time at which the expected sum-of-spins drops to a square-root of the volume, where intuitively it is absorbed within the normal fluctuations of the Ising measure.

\begin{maintheorem}\label{mainthm-trans}
For any $d\geq 2$ there exists $\beta_0=\beta_0(d)>0$ such that the following holds.
 Let $G$ be a transitive graph on $n$ vertices.
  For any fixed $0<\epsilon<1$, continuous-time Glauber dynamics for the Ising model on $G$ at inverse-temperature $0\leq\beta\leq\beta_0$ satisfies
$ \tmix(\epsilon) = \tcut \pm O_\epsilon(1)$. 

 In particular, the dynamics on a sequence of such graphs has cutoff with an $O(1)$-window around $\tcut$.
\end{maintheorem}

\section{Basic setup of information percolation and proof of Theorem~\ref{mainthm-trans}}\label{sec:proof-mainthm}
In this section we define and classify the information percolation clusters in their most basic form (commenting how this setup may be altered in more delicate situations such as $\beta$ close to $\beta_c$ on $\Z^d$), then reduce the proof of Theorem~\ref{mainthm-trans} to estimating the probability that a cluster is ``red''.

\subsection{Red, green and blue information percolation clusters}
The dynamics can be viewed as a deterministic function of $X_0$ and a random ``update sequence'' of the form $(J_1,U_1,t_1),(J_2,U_2,t_2),\ldots$, where $0<t_1<t_2<\ldots$ are the update times (the ringing of the Poisson clocks), the $J_i$'s are i.i.d.\ uniformly chosen sites (whose clocks ring), and the $U_i$'s are i.i.d.\ unit variables (random coin tosses).
According to this representation, one processes the updates sequentially: set $t_0=0$; the configuration $X_t$ for all $t\in [t_{i-1},t_i)$ ($i\geq 1$) is obtained by updating the site $J_i$ via the unit variable as follows: letting $\sigma = \sum_{v \sim J_i} X_{t_{i-1}}(v)$ denote the current sum-of-spins at the neighbors of $J_i$, if the coin toss satisfies
\begin{equation}\label{eq-plus-prob} U_i < e^{\beta\sigma }/(e^{\beta\sigma }+e^{-\beta\sigma}) = \tfrac12\left(1+\tanh(\beta\sigma)\right)\end{equation}
then the new spin at $J_i$ is chosen to be plus and otherwise it is set to minus.

Equivalently, one may evaluate this deterministic function backward in time rather than forward: sort the same update sequence $\{(J_i,U_i,t_i)\}$ such that $\tpluss>t_1>t_2>\ldots$, where $\tpluss$ is a designated time at which we wish to analyze the distribution $X_{\tpluss}$ (and argue it is either close or far from equilibrium). To construct $X_{\tpluss}$, we again process the updates sequentially, now setting $t_0=\tpluss$ and determining $X_t$ for all $t\in [t_{i+1},t_i)$ in step $i$, where the value of the spins at the neighbors of $J_i$ (determining $\sigma$ and the probability of $\pm1$ in the update as above) is evaluated recursively via the suffix of the update sequence.

Examining~\eqref{eq-plus-prob}, this backward simulation of the dynamics can be made more efficient, since even if all the neighbors of the site that is being updated are (say) minus, there is still a positive probability of $\frac12(1-\tanh(\beta d))$ for a plus update: We may therefore decide to first examine $U_i$, and if $U_i < \theta$ for
\begin{equation}\label{eq-def-theta}
  \theta = \theta_{\beta,d} := 1 - \tanh(\beta d)
\end{equation}
we will set the new spin to plus/minus as a fair coin toss irrespective of $\sigma$ (namely, to plus iff $U_i<\theta/2$);  otherwise, we will recursively compute the spins at the neighbors of $J_i$, and set the new spin to plus iff
\[ 0 \leq U_i - \theta  < \tfrac12(1+\tanh(\beta\sigma) - \theta) = \tfrac12(\tanh(\beta\sigma) + \tanh(\beta d))\,.\]
(The right-hand is $0,1-\theta$ for the extreme $\sigma=\pm d$, while for other values of $\sigma$ the rule depends on $U_i$.)

We have arrived at a branching process in the space-time slab: to recover $X_{\tpluss}(v)$ we track its lineage backward in time, beginning with the temporal edge in the space-time slab $V\times [0,\tpluss]$ between $(v,\tpluss)$ and $(v,t_i \vee 0)$, where $t_i$ is the time of the latest update to $v$. If no such update was encountered and this ``branch'' has survived to time $t=0$, it assumes the value of the initial configuration $X_0(v)$. Alternatively, an update at time $t_i$ has two possible consequences: if it features $U_i<\theta$---an {\em oblivious update}---the branch is terminated as the spin can be recovered via a fair coin toss using $U_i$; otherwise, we branch out to the neighbors of $v$, adding spatial edges between $(v,t_i)$ and $(u,t_i)$ for all $u\sim v$, and continue developing the update histories of each of them until the process dies out or reaches time 0.
This produces a graph in the space-time slab which we denote by $\sH_v$, and let $\sH_v(t)$ be its intersection with the slab $V \times \{t\}$ (viewed as a subset of $V$); further let $\sH_A = \bigcup_{v\in A} \sH_v$ and set $\sH_A(t)$ analogously.
The \emph{information percolation} clusters are the connected components of the graph consisting of $\sH_V$.

By definition, $X_{\tpluss}(A)$ is a deterministic function of the $U_i$'s corresponding to points in $\sH_A$ and of $X_0(\sH_A(0))$, the initial values at the intersection of $\sH_A$ with the slab $V\times\{0\}$; in particular, if $\sH_V(0)=\emptyset$ then $X_{\tpluss}$ is independent of $X_0$, and therefore its law is precisely the Ising distribution $\pi$ (for instance, in that scenario we could have taken $X_0\sim \pi$ and then $X_{\tpluss}\sim\pi$ by the invariance of $\pi$).
However, we note that waiting for a time $\tpluss$ large enough such that $\sH_V$ would be guaranteed not to survive to time $t=0$ with high probability is an overkill; the correct mixing time is the point at which $|\sH_V(0)|\asymp\sqrt{n}$, whence the effect of $X_0$ on $X_{\tpluss}$ would be absorbed in the normal fluctuations of $\pi$.
\begin{remark*}
The above defined rule for developing the update histories either terminated a lineage or branched it to its $d$ neighbors. In different applications of the method, it is crucial to appropriately select other rules with the correct marginal of the heat-bath dynamics.

For instance, in~\cite{LS5} we prove results {\em \`{a} la} Theorem~\ref{mainthm-trans} on any graph (including, e.g., expanders) and any $\beta < \kappa / d$ for an absolute constant $\kappa>0$---the correct dependence on $d$ up to the value of $\kappa$---by selecting $k=0,1,\ldots,d$ with probability $p_k$ and then deciding the new spin via a function of a uniformly chosen $k$-element subset of the neighbors, with the probabilities $p_k$ (and the corresponding functions to be applied) following from a discrete Fourier expansion of the original Glauber dynamics update rule.

Furthermore, in~\cite{LS4}, to analyze $\Z^d$ arbitrarily close to $\beta_c$, instead of describing such a rule explicitly one relies on the \emph{existence} of an efficient rule thanks to the exponential decay of spin-spin correlations.
\end{remark*}
\begin{example*}
In the Ising model on $\Z_n$, with probability $\theta$ we assign a uniform $\pm1$ spin, and with probability $1-\theta$ we expose $\sigma\in\{0,\pm2\}$ and select the consensus spin in case $\sigma=\pm2$ or a uniform $\pm1$ spin in case $\sigma=0$. Hence, an equivalent rule would be to terminate the branch with probability $\theta$ and otherwise to select a uniform neighbor and copy its spin, so $\sH_v$ is merely a continuous-time simple random walk killed at rate $\theta$, and $\sH_V$ consists of $n$ coalescing random walks killed at rate $\theta$.
The probability that $\sH_v(0)\neq\emptyset$ (survival from time $\tpluss$ to time $t=0$) is then $\exp(-\theta\tpluss)$, and we would have $\E|\sH_V(0)|\asymp \sqrt{n}$ at $\tpluss= (2\theta)^{-1}\log n+O(1)$, which is exactly the cutoff location.
\end{example*}

The key to the analysis is to classify the information percolation clusters to three types, where one of these classes will be revealed (conditioned on), and the other two will represent competing factors, which balance exactly at the correct point of mixing. In the basic case, the classification is as follows:
\begin{compactitem}
\item A cluster $\cC=\sH_A$ is \red\ if, given the update sequence, its final state $X_{\tpluss}(A)$ is a nontrivial function of the initial configuration $X_0$ (in particular, it survives to time $t=0$, i.e., $\sH_A(0)\neq\emptyset$).
\item A cluster $\cC=\sH_A$ is \blue\ if $A$ is a singleton ($A=\{v\} $ for some $v\in V$) and its history does not survive to time zero (in particular, $X_{\tpluss}(A)$ does not depend on $X_0$).
\item Every other cluster $\cC$ is \green.
\end{compactitem}
Let $V_\red= \{ v : \sH_v \subset\cC \in \red\}$ be the set of all vertices whose histories belong to red clusters, and let $\sH_\red = \sH_{V_\red}$ be their collective history (similarly for blue/green).
By a slight abuse of notation, we write $A\in\red$ to denote that $\sH_A\in\red$ (similarly for blue/green), yet notice the distinction between $A\in\red$ and $A\subset V_\red$ (the former means  that $\sH_A$ is a full red cluster, rather than covered by ones).

\begin{remark*}
Various other classifications to red/green/blue can be used so that \red\ captures the entire dependence on the initial configuration and \blue\ forms a product-measure. For instance, in order to study the effect of initial conditions in~\cite{LS5} we let a cluster be red if at least two branches of it survive to time $t=0$ and \emph{coalesce upon continuing to develop their history} along $t\in(-\infty,0]$; and in order to carry the analysis of~\cite{LS4} up to $\beta_c$, the classification is  more delicate, and involves the histories along a burn-in phase near time $\tpluss$ which allows one to amplify the subcritical nature of the clusters.
\end{remark*}

Note that if $\{v\}\in\blue$ then, by definition and symmetry, the distribution of $X_{\tpluss}(v)$ is uniform $\pm1$. On the other hand, while the spins of a green cluster $\cC$ are also independent of $X_0$, its spin-set at time $\tpluss$ can have a highly nontrivial distributions due to the dependencies between the intersecting update histories. It is these green clusters that embody the complicated structure of the Ising measure.

\begin{example*}
In the Ising model on $\Z_n$, as explained above, an information percolation cluster corresponds to a maximal collection of random walks that coalesced. A cluster is red if it survived to time $0$ (the rule of copying the spin at the location of the walk guarantees a nontrivial dependency on $X_0$); it is blue if the random walk started at $v$ dies out before coalescing with any other walk and before reaching time $t=0$; and it is green otherwise. Observe that the sites of a green cluster at time $\tpluss$ all have the same spin---a uniform $\pm1$ spin, independent of $X_0$---and the probability of $u,v$ belonging to the same green cluster decays exponentially in $|u-v|$ (as the walks become more likely to die out than merge).
\end{example*}

As the green clusters demonstrate the features of the Ising measure, it is tempting to analyze them in order to understand $X_{\tpluss}$. The approach we will follow does the opposite: we will \emph{condition on the entire set of histories $\sH_{V_\green}$} (or $\sH_\green$ for brevity), and study  the remaining (red and blue) clusters.

As we hinted at when saying above that if $\sH_V(0)=\emptyset$ then $X_{\tpluss}\sim \pi$, the Ising measure $\pi$ can be perfectly simulated via developing the histories backwards in time until every branch terminates without any special exception at time $t=0$. (This would be equivalent to taking $\tpluss$ larger and larger until $\sH_V$ would be guaranteed not to survive, essentially as in the ingenious \emph{Coupling From The Past} method of perfect simulation due to Propp and Wilson~\cite{PW}). Thus, we can couple $X_{\tpluss}$ to $\pi$ via the same update sequence, in which case the distribution of $V_\green$---albeit complicated---is identical (green clusters never reach $t=0$), allowing us to only consider blue/red clusters (however in a difficult conditional space where percolation clusters are forbidden from touching various parts of the space-time slab).

The second key is to keep the blue clusters---which could have been coupled to $\pi$ just like the green ones, as they too do not survive to time $t=0$---in order to water down the effect of the red clusters (which, by themselves, are starkly differently under $X_{\tpluss}$ and $\pi$). We will show that, conditioned on $\sH_\green$, the measures $\pi$ and $X_{\tpluss}$ are both very close to the uniform measure (and therefore to each other), i.e., essentially as if there were no red clusters at all and every $v\in V\setminus V_\green$ belonged to $V_\blue$.
Indeed, conditioning on the green clusters replaces the Ising measure by a contest between blue and red clusters: at large $\tpluss$, the effect of \red\ is negligible and we get roughly the uniform measure; at small $\tpluss$, the effect of \red\ dominates and $X_0$ will have a noticeable affect on $X_{\tpluss}$; the balancing point $\tcut$ has the \red\ clusters make an effect on $X_{\tpluss}$, but just within the normal fluctuations of $\pi$.

Showing that the effect of the red clusters is negligible just beyond the cutoff location will be achieved via a clever lemma of Miller and Peres~\cite{MP} that bounds the $L^2$-distance of a measure from the uniform measure in terms of a certain exponential moment. In our setting, this reduces the problem to showing:
\begin{equation}
  \label{eq-exp-moment-bound}
  \E \left[ 2^{|V_\red \cap V_\red'|} \mid \sH_\green\right] \approx 1\quad\mbox{ in probability as }n\to\infty\,,
\end{equation}
where $V_\red$ and $V'_\red$ are independent instances of the vertices whose histories are part of red clusters.

\begin{example*}
Recall the coalescing random walks representation for the information percolation clusters of the Ising model on $\Z_n$, and suppose we wanted to estimate $\E[2^{|V_\red\cap V_\red'|}]$, i.e., the left-hand of~\eqref{eq-exp-moment-bound} \emph{without} the complicated conditioning on $\sH_\green$.
Then $v\in V_\red$ iff its random walk survives to time $t=0$, which has probability $e^{-\theta \tpluss}$. By the independence, $\P(v\in V_\red\cap V_\red') = e^{-2\theta\tpluss}$. If the sites were independent (they are not of course, but the intuition is still correct), then $\E[2^{|V_\red\cap V_\red'|}]$ would break into $\prod_v \E[ 1 + \one_{\{v\in V_\red \cap V_\red'\}} ] \leq \exp(n e^{-2\theta \tpluss})$, which
for $\tpluss = (2\theta)^{-1}\log n + C$ is at most $\exp(e^{-2\theta C})$ that approaches $1$ as we increase the constant $C>0$ in $\tpluss$. (The actual calculation of the exponential moment given $\sH_\green$, especially at very low temperatures, requires quite a bit more care.)
\end{example*}

The key to obtaining the bound on the exponential moment in~\eqref{eq-exp-moment-bound}, which is the crux of the proof, is estimating a conditional probability that $A\in\red$, in which we condition not only on $\sH_\green$ but on the entire collective histories of every vertex outside of $A$, and that $A$ itself is either the full intersection of a red cluster with the top slab $V\times\{\tpluss\}$ or a collection of blue singletons.
Formally, let
\begin{equation}
  \label{eq-Psi-def}
  \Psi_A = \sup_{\sH_A^-} \P\left(A\in \red \mid \sH_A^-\,,\,\{A\in\red\} \cup \{A \subset V_\blue\}\right) \,,
\end{equation}
where
\[
\sH_A^- = \left\{\sH_v(t) : v\notin A\,,\,t\leq \tpluss\right\}
\,,\]
noting that, towards estimating the probability of $A\in\red$, the effect of conditioning on $\sH_A^-$
amounts to requiring that $\sH_A$ must not intersect $\sH_A^-$.

\begin{lemma}\label{lem-Psi}
For any $d\geq2$ and $\lambda>0$ there exist $\beta_0, C_0>0$ such that if $\beta < \beta_0$ then for any $A\subset V$ and large enough $n$, the conditional probability that $A\in\red$ at time $\tpluss$ satisfies
\[
\Psi_A \leq C_0\, \sm_{\tpluss}\, e^{- \lambda \anim(A) } \,,
\]
where $\anim(A)$ is the size of the smallest connected subgraph containing $A$.
\end{lemma}

For intuition, recall that $A\in\red$ if the histories $\{\sH_v : v\in A\}$ are all connected and $X_{\tpluss}(A)$ is a nontrivial function of $X_0$ (in particular, $\sH_A(0)\neq\emptyset$).
This is closely related to the probability that for a single $v\in A$ we have that $X_{\tpluss}(v)$ is a nontrivial function of $X_0$ (i.e., $\sH_v$ survives to time $t=0$ and creates a nontrivial dependence on $X_0$, whence the connected component of $(v,\tpluss)$ is a red cluster). Indeed, the probability of the latter event is exactly
\begin{equation}\label{eq-mag-diff}
\sm_{\tpluss} = \E X_{\tpluss}^+(v) = \P(X_{\tpluss}^+(v) \neq X_{\tpluss}^-(v) )\,,
\end{equation}
which explains the term $\sm_{\tpluss}$ in Lemma~\ref{lem-Psi}. The extra term $\exp(-\lambda \anim(A))$ is due to the requirement that the histories of $A$ must spatially connect, thus the projection of the cluster on $V$ is a connected subgraph containing $A$ (whose size is at least $\anim(A)$ by definition).

\subsection{Upper bound modulo Lemma~\ref{lem-Psi}}
Our goal is to show that $d_\tv(\tpluss)<\epsilon$ for $\tpluss = \tcut+\scut$ with a suitably large $\scut>0$, where $d_\tv(t) = \max_{x_0}\|\P_{x_0}(X_t\in\cdot) - \pi\|_\tv $.
By Jensen's inequality, $\|\psi-\phi\|_\tv \leq \E \big[\| \psi(\cdot\mid Z) - \phi(\cdot\mid Z)\|_\tv\big]$ for any two distributions $\psi$ and $\phi$ on a finite probability space $\Omega$ and random variable $Z$.
Applied with $\sH_\green$ playing the role of $Z$, and letting $X'_0\sim \pi$,
\begin{align*}
 d_\tv(t) &\leq  \max_{x_0} \E\Big[\left\|\P_{x_0}(X_t\in \cdot \mid \sH_\green) - \P_{X'_0}(X_t\in \cdot \mid \sH_\green)\right\|_\tv \Big]\\
 &\leq \sup_{\sH_\green} \max_{x_0} \left\|\P_{x_0}(X_t\in \cdot \mid \sH_\green) - \P_{X'_0}(X_t\in \cdot \mid \sH_\green)\right\|_\tv \,.
 \end{align*}
As explained above, since $X_t(V_\green)$ is independent of $X_0$ we can couple it with the chain started at $X'_0 \sim \pi$, whence the projection onto $V\setminus V_\green$ does not decrease the total-variation distance, and so
 \begin{align}
d_\tv(t) &\leq \sup_{\sH_\green} \max_{x_0} \left\|\P_{x_0}(X_t(V \setminus V_\green) \in \cdot \mid \sH_\green) - \P_{X'_0}(X_t(V \setminus V_\green) \in \cdot \mid \sH_\green)\right\|_\tv \nonumber\\
&\leq 2 \sup_{\sH_\green} \max_{x_0} \left\|\P_{x_0}(X_t(V \setminus V_\green) \in \cdot \mid \sH_\green) - \nu_{V\setminus V_\green}\right\|_\tv\,,\label{eq-exp-moment-1}
\end{align}
where $\nu_A$ is the uniform measure on configurations on the sites in $A$.

We now appeal to the aforementioned  lemma of Miller and Peres~\cite{MP} that shows that, if a measure $\mu$ on $\{\pm1\}^V$ is given by sampling a variable $R\subset V$ and using an arbitrary law for its spins and a product of Bernoulli($\frac12$) for $V\setminus R$, then the $L^2$-distance of $\mu$ from the uniform measure is at most $\E 2^{|R\cap R'|}-1$ for i.i.d.\ copies $R,R'$. (See~\cite{LS4}*{Lemma~4.3} for a generalization of this to a product of general measures that is imperative for the information percolation analysis on $\Z^d$ at $\beta$ near $\beta_c$.)
\begin{lemma}[\cite{MP}]\label{lem:MP}
Let $\Omega=\{\pm1\}^V$ for a finite set $V$. For each $R\subset V$, let $\phi_R$ be a measure on $\{\pm1\}^{R}$.
Let $\nu$ be the uniform measure on $\Omega$, and let $\mu$ be the measure on $\Omega$ obtained by sampling a subset $R\subset V$ via some measure $\tilde{\mu}$, generating the spins of
 $R$ via $\phi_R$, and finally sampling $V \setminus R$ uniformly. Then
\[ \left\|\mu - \nu\right\|^2_{L^2(\nu)}   \leq
\E2^{\left|R\cap R'\right|}-1\qquad\mbox{ where $R,R'$ are i.i.d.\ with law $\tilde{\mu}$}\,.\]
\end{lemma}

Applying this lemma to the right-hand side of~\eqref{eq-exp-moment-1},
while recalling that any two measures $\mu$ and $\nu$ on a finite probability space satisfy $\|\mu-\nu\|_\tv  =\frac12 \|\mu-\nu\|_{L^1(\nu)} \leq \frac12 \|\mu-\nu\|_{L^2(\nu)}$,
we find that
\begin{align}\label{eq-exp-moment-2}
d_\tv(\tpluss) \leq \Big( \sup_{\sH_\green} \E
\left[2^{\left|V_\red \cap V_{\red'}\right|} \;\big|\; \sH_\green\right] - 1 \Big)^{1/2} \,,
\end{align}
where $V_\red$ and $V_{\red'}$ are i.i.d.\ copies of the variable $\bigcup\{v \in V : \sH_v \in \red\}$. We will reduce the quantity $|V_{\red}\cap V_{\red'}|$ to one that involves the $\Psi_A$ variables defined in~\eqref{eq-Psi-def} (which will thereafter be controlled via Lemma~\ref{lem-Psi}) using the next lemma, whose proof is deferred to~\S\ref{sec-coupling} below.

\begin{lemma}\label{lem-coupling}
Let $\{ Y_{A} : A\subset V\}$ be a family of independent indicators satisfying $\P(Y_{A}=1) = \Psi_{A}$.
The conditional distribution of $V_\red$ given $\sH_\green$ can be coupled such that
\[ \{A : A \in \red \} \subset \{A : Y_A = 1 \} \,.\]
\end{lemma}
The following corollary is then straightforward (see \S\ref{sec-coupling} for its short proof).
\begin{corollary}\label{cor-coupling2}
Let $\{ Y_{A,A'} : A,A'\subset V\}$ be a family of independent indicators satisfying
\begin{equation}
  \label{eq-YAA'-def}
  \P(Y_{A,A'}=1) = \Psi_{A}\Psi_{A'}\quad\mbox{ for any $A,A'\subset V$}\,.
\end{equation}
The conditional distribution of $(V_\red,V_{\red'})$ given $\sH_\green$ can be coupled to the $Y_{A,A'}$'s such that
\[ \left|V_\red \cap V_{\red'}\right| \preceq \sum_{A\cap A'\neq \emptyset} |A\cup A'| Y_{A,A'}\,.\]
\end{corollary}
Relaxing $|A\cup A'|$ into $|A|+|A'|$, we get \begin{align*}
\sup_{\sH_\green}\E\left[2^{|V_\red \cap V_{\red'}|} \;\big|\; \sH_\green\right] &\leq  \E\left[ 2^{\sum_{A\cap A'\neq\emptyset} (|A|+|A'|) Y_{A,A'}}\right] = \prod_{A\cap A'\neq \emptyset} \E\left[ 2^{(|A|+|A'|)Y_{A,A'}}\right]\,,
\end{align*}
with the equality due to the independence of the $Y_{A,A'}$'s. By the definition of these indicators in~\eqref{eq-YAA'-def},
\begin{align*}
\prod_{A\cap A'\neq \emptyset} \E\left[ 2^{(|A|+|A'|)Y_{A,A'}}\right] \leq \prod_v \!\prod_{\substack{A,A'\\ v\in A\cap A'}} \!\left(\big(2^{|A|+|A'|}-1\big) \Psi_{A}\Psi_{A'}+1\right)
\leq e^{n \left(\sum_{A \ni v} 2^{|A|} \Psi_{A}\right)^2}\,,
\end{align*}
and so, revisiting~\eqref{eq-exp-moment-2}, we conclude that
\begin{equation}\label{eq-req-exp-bound}
 d_\tv(\tpluss) \leq  \left( e^{n\left(\sum_{A \ni v} 2^{|A|} \Psi_{A}\right)^2} - 1\right)^{1/2} \wedge 1 \leq \sqrt{2 n} \sum_{A \ni v} 2^{|A|} \Psi_{A}\,,
 \end{equation}
where we used that $e^x-1\leq 2x$ for $x\in[0,1]$. Using Lemma~\ref{lem-Psi} with
$ \lambda= \log(4ed) $
we find that
\begin{align*}
\sum_{A \ni v} 2^{|A|} \Psi_A \leq C_0 \sm_{\tpluss}  \sum_k \sum_{\substack{A \ni v \\ \anim(A)=k}} 2^k e^{-\lambda k}
\leq C_0\sm_{\tpluss} \sum_k (2 e d  e^{-\lambda})^k \leq 2 C_0 \sm_{\tpluss}\,,
\end{align*}
for some $C_0=C_0(d)>0$, and going back to~\eqref{eq-req-exp-bound} shows that
$
 d_\tv(\tpluss) \leq 2\sqrt{2} C_0\, \sm_{\tpluss}\sqrt{n}
$.

The upper bound is concluded by the submultiplicativity of the magnetization (see~\cite{LS5}*{Claim~3.3}):
\begin{equation}\label{eq-MagnetDecay}
e^{-t} \sm_{t_0}\leq
\sm_{t_0+t} \leq e^{-(1-\beta d) t }\sm_{t_0}  \qquad\mbox{ for any $t_0,t\geq 0$}
\,,\end{equation}
since $\sm_{t_\sm} = 1/\sqrt{n}$ (recall~\eqref{eq-t*-def}) and so taking $\tpluss = t_\sm + C$ for a suitable $C(\epsilon)>0$ yields $d_\tv(\tpluss)<\epsilon$.
\qed

\subsection{Proof of Lemma~\ref{lem-coupling}}\label{sec-coupling}
We claim that given $\sH_\green$, if one arbitrarily orders all distinct subsets $A\subset V\setminus V_\green$ as $\{A_i\}_{i\geq 1}$ and then successively exposes whether $\{A_i\in\red\}$, denoting the associated filtration by $\cF_i$, then $\P(A_i\in\red \mid\cF_{i-1})\leq\Psi_{A_i}$. To see this, first increase $\P(A_i\in\red\mid\cF_{i-1})$ into $\P(A_i\in\red\mid \{A_i \in \red\}\cup\{A_i\subset V_\blue\},\cF_{i-1})$, and then further condition on a worst-case $\sH_{A_i}^-$.
The latter subsumes any information in $\one_{\{A_j\in\red\}}$ for $A_j$ disjoint from $A_i$; the former means that $A_i$ is either the full intersection of a red cluster with $V\times\{\tpluss\}$ or a collection of blue singletons, either way implying that any $A_j$ intersecting $A_i$ must satisfy $A_j\notin\red$. Altogether, the events $\{A_j\in\red : j<i\}$ are measurable under the conditioning on $\{A_i\in\red\}\cup\{A_i\subset V_\blue\},\sH_{A_i}^-$, and we arrive at $\Psi_{A_i}$ which immediately implies the stochastic domination.
\qed

\begin{proof}[\textbf{\em Proof of Corollary~\ref{cor-coupling2}}]
By Lemma~\ref{lem-coupling} we can stochastically dominate $\{A\in\red\}$ and $\{A'\in\red'\}$ by two independent sets of indicators $\{Y_A\}$ and $\{Y'_A\}$.
Let $\{(A_l,A'_l)\}_{l\geq 1}$ denote all pairs of intersecting subsets ($A,A'\subset V \setminus V_\green$ with $A\cap A'\neq\emptyset $) arbitrarily ordered, and associate each pair with a variable $R_l$ initially set to 0. Process these in order: If we have not yet found $R_j = 1$ for some $j<l$ with either $A_j\cap A_l\neq\emptyset $ or $A'_j\cap A'_l\neq\emptyset$, then set $R_l = \one_{\{A_l\in\red,\, A'_l\in\red'\}}$ (otherwise skip this pair, keeping $R_l=0$).

Noting $\P(R_l = 1 \mid \cF_{l-1}) \leq \P(Y_{A_l}=1,Y'_{A'_l} = 1) = \Psi_{A_l} \Psi_{A'_{l}}$ (as testing $R_l = 1$ means we received only negative information on $\{Y_{A_l}=1\}$ and $\{Y'_{A'_l}=1\}$) gives the coupling $\{l:R_l=1\} \subset \{l: Y_{A_l,A'_l}=1\}$.

Completing the proof is the fact that if $v\in V_\red\cap V_{\red'}$ then there are subsets $A_l,A'_l$ containing it such that $Y_{A_l}$ and $Y'_{A'_l}$, in which case either $R_l=1$  or we will have an earlier $R_j=1$ for a pair involving $A_j=A_l$ or $A_j'=A_l'$ (nontrivial intersections with $A_l$ or $A'_l$ will not be red), whence $v\in A_j\cup A'_j$.
\end{proof}

\subsection{Lower Bound}
Recalling the choice of $\tcut$ as the point such that $\sm_{\tcut} = 1/\sqrt{n}$,
we let the sum of spins $f(\sigma) = \sum_v \sigma(v)$ be our distinguishing statistic at time $\tminuss = \tcut-\scut$.
  Putting $Y = f\big(X^+_{\tminuss}\big)$ for the dynamics $(X_t^+)$ started from all-plus, by~\eqref{eq-MagnetDecay} we have
\begin{align}
  \label{eq-E[Y]}
\E Y = n \sm_{\tminuss} \geq e^{2(1-\beta d) \scut} \sm_{\tcut} \geq e^{\scut} \sqrt{n}
\end{align}
(the last inequality using $\beta d \leq \frac12$).
For the variance estimate we use the fact (see~\cite{LS5}*{Claim~3.4}
) that for some constants $\beta_0=\beta_0(d)>0$ and $\gamma=\gamma(d)>0$, if $\beta < \beta_0$ then
\begin{equation}\label{eq:correlation}
\sum_u \Cov(X_{t}(u), X_t(v)) \leq \gamma \quad\mbox{ for any $X_0$, $t>0$ and $v\in V$}\,.
\end{equation}
From this inequality it follows that $ \var\left(Y\right) \leq \gamma n$, and in light of~\eqref{eq-E[Y]}, $\P( Y \geq \E Y / 2) > 1-\epsilon/2$ by Chebyshev's inequality provided $\scut=\scut(\epsilon)$ is chosen large enough.

On the other hand, if $X'\sim \pi$ then $\E [f(X')] = 0$ (as $\E[X'(v)]=0$ for any $v$), while $\var (f(X')) \leq \gamma n$  by the decay of correlation of the Ising measure. By Chebyshev's inequality, for any large enough $\scut=\scut(\epsilon)$ we have $\P(f(X') \geq \E Y / 2) < \epsilon/2$. Altogether, $d_\tv(\tminuss) \geq 1-\epsilon$,  as required.
\qed

\section{Analysis of the information percolation clusters}\label{sec:cluster-analysis}

The delicate part of the bounding $\Psi_A$ is of course the conditioning: since $A$ is either the interface of a complete red cluster, or a collection of blue singletons, its combined histories must avoid $\sH_A^-$. This immediately implies, for instance, that if any branches of $\sH_A^-$ should extend to some point $(v,t)\in A\times [0,\tpluss)$ in the space-time slab, necessarily the branch of $v$ must receive an update along the interval $(t,\tpluss]$ to facilitate avoiding that branch.
Our concern will be such a scenario with $t\in (\tpluss-1,\tpluss]$, since for values of $t$ extremely close to $\tpluss$ this event might be extremely unlikely (potentially having probability smaller than $\exp[-O(|A|)]$, which would we not be able to afford). For a subset $A'\subset A$ and a sequence of times $\{s_u\}_{u\in A'}$ with $s_u\in(\tpluss-1,\tpluss]$, define
\begin{equation}
  \label{eq-U-def}
  \cU = \cU(A',\{s_u\}_{u\in A'}) = \bigcap_{u\in A'}\big\{ \mbox{$u$ receives an update along $(s_u,\tpluss]$}\big\}\,.
\end{equation}
We will reduce the conditioning on $\sH_A^-$ to an appropriate event $\cU$, and thereafter we would want to control $\sH_A$, the collection of all histories from $A$, both spatially, measured by its branching edges
\[\chi(\sH_A) = \#\left\{\big((u,t),(v,t)\big)\in\sH_A\;:\;u v \in E(G)\,,\,t\in(0,\tpluss]\right\}\]
(i.e., those edges that correspond to a site branching out to its neighbors via a non-oblivious update), and temporally, as measured by the following length quantity:
\begin{equation}\label{eq-len-def}
 \len(\sH_A) = \sum_{u\in V} \int_0^{\tpluss} \one_{\{(u,t)\in\sH_A\}}dt\,.
\end{equation}
 The following lemma bounds an exponential moment of $\chi(\sH_A)$ and $\len(\sH_A)$ under any conditioning on an event $\cU$ as above, and will be central to the proof of Lemma~\ref{lem-Psi}.
\begin{lemma}
  \label{lem-todo}
For any $d\geq 2$, $0<\eta<1$ and $\lambda>0$ there are $\beta_0(d,\eta,\lambda)>0$ and $\alpha(d,\eta)>0$ such that the following holds. If $\beta<\beta_0$ then for any subset $A$,
  \[ \sup_{\cU} \E\left[ \exp\left( \eta \len(\sH_A) +  \lambda \,\chi(\sH_A)\right) \given \cU\right] <  \exp\left(\alpha |A|\right)\,,\]
  where the event $\cU$ is as defined in~\eqref{eq-U-def}.
\end{lemma}

\subsection{Proof of Lemma~\ref{lem-Psi}}
For a given subset $S\subset V$, let $\red^*_S$ denote the red clusters that arise when exposing the joint histories of $\sH_S$ (as opposed to all the histories $\sH_V$),
noting the events $\{A\in\red\}$ and $\{A\in\red_A^*\}\cap\{\sH_A\cap\sH_A^-=\emptyset\} $ are identical
(so that $A$ would be the interface of a full red cluster).
Similarly define $\blue^*_S$, and by the same reasoning $\{A\subset V_\blue\} = \{A\subset V_{\blue^*_A}\}\cap\{\sH_A\cap\sH_A^-=\emptyset\}$.

Next, given $\sH_A^-=\cX$, let $ s_u = s_u(\cX) = \max\left\{ s \;:\; (u,s) \in \cX\right\}$ be the latest most time at which $\cX$ contains $u\in A$, and recall from the discussion above that any $u$ with $s_u\leq\tpluss$ must receive an update along $(s_u,\tpluss]$ in order to avoid $\cX$.  Thus, writing $ A' = \left\{ u \in A \;:\; s_u > \tpluss-1\right\}$
and defining $\cU(A',\{s_u\}_{u\in A'})$ as in~\eqref{eq-U-def}, we find $\P\left(A\in\red\mid\sH_A^-=\cX,\{A\in\red\}\cup\{A\subset V_\blue\}\right)$ to be equal to
\[\frac{\P(A\in\red_A^*\,,\,\sH_A\cap\cX=\emptyset\,,\,\cU \mid \sH_A^-=\cX)}{\P(\{A\in\red^*_A\}\cup\{A\subset V_{\blue^*_A}\}\,,\,\sH_A\cap\cX=\emptyset\,,\,\cU \mid \sH_A^-=\cX)}\,,\]
which, since both of the events $A\in\red_A^*$ and $A\subset V_{\blue^*_A}$ are $\sH_A$-measurable, equals
\[\frac{\P(A\in\red_A^*\,,\,\sH_A\cap\cX=\emptyset\mid \cU)}{\P(\{A\in\red^*_A\}\cup\{A\subset V_{\blue^*_A}\}\,,\, \sH_A\cap\cX=\emptyset\mid \cU)}\,.\]
The numerator is at most $\P(A\in\red_A^*\mid\cU)$. As for the denominator, it is at least the probability that, in the space conditioned on $\cU$, every $u\in A$ gets updated in the interval $(s_u\vee \tpluss-1,\tpluss]$ and the last such update (i.e., the first we expose when revealing $\sH_u$) is oblivious (implying $A\subset V_{\blue_A^*}$)---which is $\theta^{|A|}(1-1/e)^{|A\setminus A'|}$. As this is at least $e^{-|A|}$ for small enough $\beta$ (recall the definition of $\theta$ in~\eqref{eq-def-theta}),
\begin{align}
  \label{eq-psi-bound-given-U}
\Psi_A \leq e^{|A|} \P(A\in\red_A^* \mid \cU)\,.
\end{align}
To estimate $\P(A\in\red_A^*\mid\cU)$, let us expose the history of $A$ backwards from $\tpluss$ along the first unit interval $(\tpluss-1,\tpluss]$ (moving us past the information embedded in $\cU$), then further on until reaching time $T$ where $\sH_A(T)$ coalesces to a single point or we reach time $T=0$. For $A\in\red_A^*$ to occur, in either situation the projections of $X^+_T$ and $X^-_T$ on the subset $\sH_A(T)$ must differ (otherwise $X_{\tpluss}(A)$ will not depend nontrivially on $X_0$). If $T=0$ this trivially holds, and if $T>0$ then, as $\sH_A(T)=\{(w,T)\}$ for some vertex $w$ and we did not expose any information on the space-time slab $V\times[0,T]$, the probability of this event is exactly $\sm_T$.
Furthermore, using~\eqref{eq-MagnetDecay} we have $\sm_{T} \leq e^{\tpluss-T} \sm_{\tpluss}$, and by definition, along the interval $[T,\tpluss-1]$ there are at least two branches in $\sH_A$, so $\len(\sH_A) \geq 2(\tpluss-1-T)$; thus, $\sm_T \leq e^{\frac12\len(\sH_A)+1}\sm_{\tpluss}$.  Also, on the event $A\in\red_A^*$ the histories $\sH_A$ must all join by time $T$, and thus $\chi(\sH_A) \geq \anim(A) - 1$ since $\sH_A$ must observe at least that many branching edges for connectivity.

In conclusion, for any $v\in A$,
\begin{align*}
\P(A\in\red_A^* \mid \cU) &\leq   e \,\sm_{\tpluss} \E\left[ \one_{\{\chi(\sH_A) \geq \anim(A) - 1\}} e^{\frac12\len(\sH_A)} \mid \cU\right]\,,
\end{align*}
which for any $\lambda>0$ and $\alpha>0$ is at most
\begin{align*}
 e^{1-(\lambda+\alpha+1) (\anim(A)-1)} \sm_{\tpluss} \E\left[ e^{(\lambda+\alpha+1)\chi(\sH_A) + \frac12\len(\sH_A)} \mid \cU\right]\,.
\end{align*}
Plugging in $\alpha$ as given by Lemma~\ref{lem-todo} (which we recall does not depend on the pre-factor of $\chi(\sH_A)$ in that lemma), the exponential moment above is at most $e^{\alpha |A|}$, and revisiting~\eqref{eq-psi-bound-given-U} we conclude that
\[\Psi_A \leq  \sm_{\tpluss} e^{1+(1+\alpha)|A|-(\lambda+\alpha+1)(\anim(A)-1)} \leq
 e^{\lambda+\alpha+2}\,\sm_{\tpluss} e^{-\lambda \anim(A)} \,.
\tag*{\qed}
\]

\subsection{Exponential decay of cluster sizes: Proof of Lemma~\ref{lem-todo}}
We develop the history of a set $A$ backward in time from  $\tpluss$ by exposing the space-time slab.
Let $W_s = |\sH_A(\tpluss - s)|$ count the number of vertices in the history of $A$ at time $\tpluss - s$, let $Y_s = \chi(\sH_A \cap V \times [\tpluss - s, \tpluss])$ be the total number of vertices observed by the history by time $\tpluss - s$ and let $Z_s = \sum_{u\in V} \int_{\tpluss - s}^{\tpluss} \one_{\{(u,t)\in\sH_A\}}dt$ be the total length of the history in  the time interval $[\tpluss - s, \tpluss]$ of the space-time slab.  Initially we have $(W_0,Y_0,Z_0)=(|A|,0,0)$.

Recall that the probability that an update of a vertex $v$ will branch out to its $d$ neighbors is $1-\theta$ and that with probability $\theta$ it is oblivious which observes no neighbours.  Thus we can stochastically dominate $(W_s,Y_s,Z_s)$ by a process $(\bar{W}_s,\bar{Y}_s,\bar{Z}_s)$ defined as follows. Initially, $(\bar{W}_0,\bar{Y}_0,\bar{Z}_0)=(|A|,0,0)$ and at rate $\theta\bar{W}_s$ we decrease $\bar{W}_s$ by 1 and at rate $(1-\theta)\bar{W}_s$ both $\bar{W}_s$ and $\bar{Y}_s$ increase by $d$.  The length grows as $d\bar{Z}_s = \bar{W}_s ds$.
Now consider the process,
\[
Q_s = \exp\left( \eta \bar{Z}_s + \lambda  \bar{Y}_s + \alpha \bar{W}_s \right).
\]
for some $\alpha > - \log(1-\eta)$ which does not depend on $\lambda$.  We have that
\[
\frac{d}{ds} \E[Q_s \mid Q_{s_0}] \,\Big|_{s=s_0} =  \left( \eta  + \theta (e^{-\alpha}-1) + (1-\theta) (e^{(\lambda+\alpha)d}-1) \right)\bar{W}_{s_0} Q_{s_0}
\]
which is negative provided $\theta$ is sufficiently close to 1 (guaranteed by taking $\beta_0$ sufficiently small). Then, letting $\tau$ be the first time that $\bar{W}_\tau=0$, optional stopping for the supermartingale $Q_s$ yields
\[
\E \exp\left( \eta \bar{Z}_S + \lambda  \bar{Y}_S\right) \leq \E Q_0 = \exp(\alpha |A|)\,.
\]
By the stochastic domination we have that $\E\left[ \exp\left( \eta \len(\sH_A) + \lambda \, \chi(\sH_A)\right)\right] \leq \exp(\alpha |A|)$.
Under this coupling the effect of conditioning on $\cU$ is simply to expedite updates and hence reduce the length of the process, thus the conditioning can only reduce the expectation.
\qed

\subsection*{Acknowledgements}  AS was supported by an Alfred Sloan Fellowship and NSF grants DMS-1208338,  DMS-1352013. The authors thank the Theory Group of Microsoft Research, Redmond, where much of the research has been carried out, as well as P.\ Diaconis and Y.\ Peres for insightful discussions.

\end{document}